# Spectrum of analytic continuation

D V Ingerman

**Abstract.** I will show that operator of analytic (harmonic) continuation on a lattice graph has a positive spectrum. I use a theorem about positivity of eigenvalues of totally positive matrices. I conjecture that by approximation the similar result holds in continuous case on a plane.

**1. Harmonic functions**

Let $\gamma$ be a positive function (conductivity) on the edges of a graph. We call two vertices of a graph neighbors if they are connected by an edge. A vector $u$ on the vertices of the graph is called *harmonic* or *$\gamma$-harmonic* if its value at a vertex is the weighted average of its values at the neighbors.

$$u_i = \frac{\sum_j u_j \gamma_{ij}}{\sum_j \gamma_{ij}}$$

or

$$\sum_j \gamma_{ij}(u_i - u_j) = 0$$

(analog of the Laplace-Beltrami equation

$$u_{xx} + u_{yy} = 0$$

or

$$(\gamma u_x)_x + (\gamma u_y)_y = 0$$

.)

**2. Harmonic continuation**

If the values of a harmonic function $u$ are specified at some vertices of a graph

$$a_1, a_2, ..., a_n$$

(or a region in continuous case) they may imply the values of $u$ at some other vertices (or a region)

$$b_1, b_2, ..., b_m$$

by *harmonic continuation*.
The map of harmonic continuation

$$H : u(a_i) \to u(b_j)$$

is linear.
We will show that the spectrum of $H$ for $n = m$ for the following graph is positive and conjecture that the same result holds for arbitrary planar graphs and continuous case on a plane.

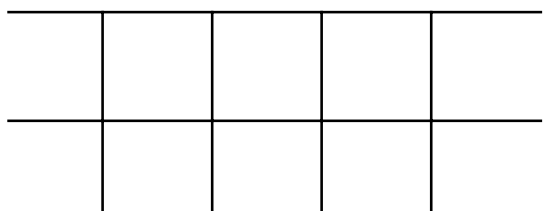

**Figure 1**. We consider harmonic functions on this lattice.

### 3. Spectrum of harmonic continuation
We will prove that the eigenvalues of *modified H* are positive by constructing a basis in which the matrix representing the linear operator of harmonic continuation *H* is *totally non-negative* (the determinants of all its square submatrices are non-negative). By [8, 10] this will imply the result.
Let modified *H* be the map from differences between values of *u* at the edges *a* to the differences between values of *u* at the edges *b*.

$$H : \nabla u(a) \to \nabla u(b)$$

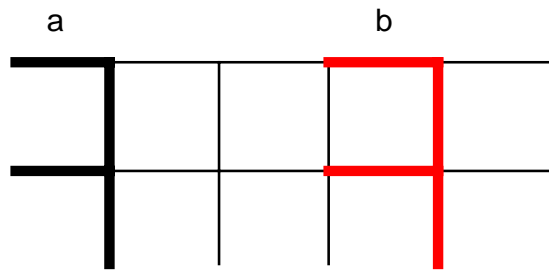

**Figure 2**. We consider harmonic continuation from black edges a to red edges b.

Lets take the differences $\nabla u$ of *u* in the directions of the arrows shown on the figure 3.

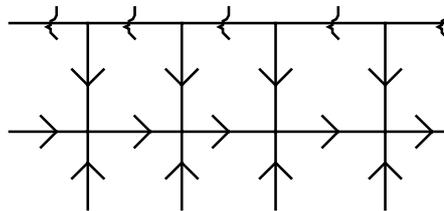

**Figure 3**. Arrows making harmonic continuation totally positive.

The directions are chosen in such a way that the harmonicity equation for *u*

$$\sum_j \gamma_{ij}(u_i - u_j) = 0$$

implies that $\nabla u(red)$ is a linear combination with **positive** coefficients of $\nabla u(black)$ in the figure 4

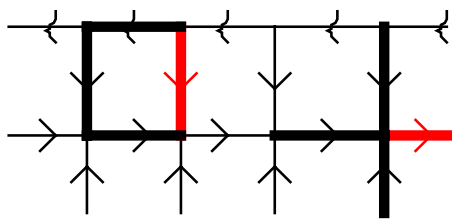

**Figure 4**. $\nabla u(red)$ is a linear combination with **positive** coefficients of $\nabla u(black)$.

This is a step in the harmonic continuation.

The dependence is obviously given by the matrix of the following form

$$\nabla u(red) = \begin{pmatrix} 1 & 0 & 0 & 0 & 0 \\ 0 & 1 & 0 & 0 & 0 \\ 0 & + & + & + & 0 \\ 0 & 0 & 0 & 1 & 0 \\ 0 & 0 & 0 & 0 & 1 \end{pmatrix} \nabla u(black)$$

which is obviously totally non-negative. The harmonic continuation $H$ from edges $a$ to edges $b$ (figure 2) is a sequence of steps represented by the matrices of the form above. The product of non-singular totally non-negative matrices is a non-singular totally non-negative matrix by Cauchy-Binet formula for determinants. Therefore $H$ as the product of the matrices of the form above is represented by a non-singular totally non-negative matrix. The spectrum of totally non-negative matrix is non-negative [8,10], therefore

$$\nabla u(b) = \lambda \nabla u(a)$$

has 5 positive solutions

$$\lambda > 0.$$

## 4. Conjecture for the continuous case

Let $\gamma$ be a positive function (conductivity) on a region $\Omega$ in a plane. A function $u$ on the region is called *harmonic* or *$\gamma$-harmonic* on $\Omega$ if it satisfies Laplace-Beltrami equation

$$(\gamma u_x)_x + (\gamma u_y)_y = 0$$

or Cauchy-Riemann system

$$\begin{cases} \gamma u_x = v_y \\ \gamma u_y = -v_x \end{cases}$$

for some function $v$.

**Conjecture 4.1**

Let $\Omega$ be a region in a plane and $\gamma$ be a positive function on $\Omega$.
Suppose
- funtion $u$ is $\gamma$-harmonic on $\Omega$
- regions $A$ and $A+shift$ belong to $\Omega$
- the values if $u$ at the region $A$ determine the vales of $u$ at the region $A+shift$ by $\gamma$-harmonic continuation

If for some $\lambda$

$$\nabla u(z + shift) = \lambda \nabla u(z)$$
$$\forall z \in A$$

then

$$\lambda > 0.$$

The proof should follow from the approximation of the continuous equations by the equations on a lattice.

## 5. Applications

The eigenvalues of the harmonic continuation can be directly measured from the Dirichlet-to-Neumann map of the graph implying simpler equations for conductivity in terms of the Dirichlet-to-Neumann map for the inverse problems.

**References**

[1] The Dirichlet to Neumann Map for a Resistor Network Curtis E B, Morrow J A *SIAM Journal on Applied Mathematics*, Vol. **51** No. 4 (Aug., 1991) pp. 1011-1029

[2] Loop-erased walks and total positivity Fomin S *Transactions of the AMS* **353** (2001) 3563-3583

[3] Optimal finite-difference grids and rational approximations of square root. I. Elliptic problems Druskin V, Ingerman D V, Knizhnerman L *Comm. Pure. Appl. Math.* **53** (2000) no. 8 pp. 1039-1066

[4] Discrete and Continuous Dirichlet-to-Neumann Maps in the Layered Case Ingerman D V *SIMA* Vol. **31** Number 6 pp. 121-123 2000

[5] Circular planar graphs and resistor networks Curtis E B, Ingerman D V, Morrow J A *Linear Algebra and Its Applications,* **283** (1998) 115-150

[6] On a Characterization of the Kernel of the Dirichlet-to-Neumann Map for a Planar Region Ingerman D V, Morrow J A *SIMA* Vol. **29** Number 1 pp. 106-115 1998

[7] Discrete and continuous inverse boundary problems on a disc Ingerman D V *Ph.D. Thesis* University of Washington 1997

[8] *Total Positivity* Karlin S Stanford University Press 1968

[9] A Comparison of Three Algorithms for the Inverse Conductivity Problem Landrum J University of Washington *REU on inverse problems* 1990

[10] *Oscillation Matrices and Kernels and Small Vibrations of Mechanical Systems* Gantmakher F R, Krein M G Moscow 1950